\documentclass{article}

\usepackage{amsmath,amssymb}
\usepackage{graphicx}
\usepackage[colorlinks=true,citecolor=black,linkcolor=black,urlcolor=blue]{hyperref}
\usepackage{url}
\usepackage{cleveref}

\def\paren#1{\left( #1 \right)}
\def\acc#1{\left\{ #1 \right\}}

\renewcommand{\le}{\leqslant}
\renewcommand{\ge}{\geqslant}

\newcommand{\arxiv}[1]{\href{http://arxiv.org/abs/#1}{\texttt{arXiv:#1}}}
\newtheorem{theorem}{Theorem}

\newtheorem{remark}[theorem]{Remark}

\title{A family of formulas with reversal of arbitrarily high avoidability index}

\author{Pascal Ochem\thanks{The author was partially supported by the ANR project CoCoGro (ANR-16-CE40-0005).}\\
\small LIRMM, CNRS, Universit\'e de Montpellier\\[-0.8ex] 
\small France\\
\small\tt ochem@lirmm.fr}
\begin{document}

\maketitle
\setcounter{footnote}{0}
\begin{abstract}
We present a family of avoidable formulas with reversal whose avoidability index is unbounded.
We also complete the determination of the avoidability index of the formulas with reversal in the 3-avoidance basis.
\end{abstract}

\section{Introduction}\label{sec:intro}
The notion of formula with reversal~\cite{CMR2017,CMR2018} is an extension of the notion of classical formula such that a variable $x$ can appear both as $x$ and $x^R$
with the convention that in an occurrence $h$ of the formula, $h(x^R)$ is the reverse (i.e., mirror image) of $h(x)$.
The avoidability index $\lambda(F)$ of a formula with reversal $F$ is the minimum number of letters contained in an infinite word avoiding $F$.

Currie, Mol, and Rampersad~\cite{CMR2017} have asked if there exist formulas with reversal with arbitrarily large avoidability index.
They considered the formula $\psi_k=xy_1y_2\ldots y_kx\cdot y_1^R\cdot y_2^R\cdot\ldots\cdot y_k^R$ and
obtained that $\lambda(\psi_1)=4$, $\lambda(\psi_2)=\lambda(\psi_3)=\lambda(\psi_6)=5$, $5\le\lambda(\psi_4)\le6$, $5\le\lambda(\psi_5)\le7$,
$4\le\lambda(\psi_k)\le6$ if $k\ge7$ and $k\not\equiv0\pmod3$, and $4\le\lambda(\psi_k)\le5$ if $k\ge9$ and $k\equiv0\pmod3$.
They conjecture that $\lambda(\psi_k)=5$ for all $k\ge2$.
Computational experiments suggest that the upper bound $\lambda(\psi_k)\le5$ for $k\ge3$ is witnessed
by the image of every $(7/4^+)$-free ternary word by the following $(k+3)$-uniform morphism where $k=3t+i$, $t\ge1$, and $0\le i\le2$.
$$
\begin{array}{ll}
\texttt{0}\to&(\texttt{012})^{t+1-i}(\texttt{0123})^i\\
\texttt{1}\to&(\texttt{013})^{t+1-i}(\texttt{0134})^i\\
\texttt{2}\to&(\texttt{014})^{t+1-i}(\texttt{0142})^i\\
\end{array}
$$
We give a positive answer to their original question with Theorem~\ref{phi} below. 

We define the formula $\phi_k=x_0x_1\cdot x_1x_2\cdot\ldots\cdot x_{k-1}x_0\cdot x_0^R\cdot x_1^R\cdot\ldots\cdot x_{k-1}^R$.

\begin{theorem}\label{phi}
For every fixed $b$, there exists $k$ such that $b<\lambda(\phi_k)\le k+1$.
\end{theorem}

Currie, Mol, and Rampersad~\cite{CMR2018} have also determined the 3-avoidance basis for formulas with reversal,
which contains the minimally avoidable formulas with reversal on 3 variables.
They obtained several bounds on the avoidability index of the formulas with reversal in the 3-avoidance basis.
The next two results finish the determination of the avoidability index of these formulas.

\begin{theorem}\label{xyzyx}
The formulas $xyzyx\cdot zyxy^Rz$, $xyzyx\cdot zy^Rxyz$, $xyzyx\cdot zy^Rxy^Rz$, $xyzy^Rx\cdot zyxy^Rz$, and $xyzy^Rx\cdot zy^Rxyz$ are simultaneously 2-avoidable.
\end{theorem}

\begin{theorem}\label{xyzx}
The formulas $xyzx\cdot yzxy\cdot z^R$ and $xyzx\cdot yz^Rxy$ are simultaneously 3-avoidable.
\end{theorem}

\Cref{phi,xyzyx,xyzx} are proved in \Cref{pphi,pxyzyx,pxyzx}.

A word $w$ is \emph{$d$-directed} if for every factor $f$ of $w$ of length $d$, the word $f^R$ is not a factor of $w$.
\begin{remark}\label{remark}
If a $d$-directed word contains an occurrence $h$ of $x\cdot x^R$, then $|h(x)|\le d-1$.
\end{remark}

In order to express the simultaneous avoidance of similar formulas, as in \Cref{xyzyx,xyzx}, we introduce the notation $x^U$
to represent equality up to mirror image. That is, if $h(x)=w$, then $h(x^R)=w^R$ and $h(x^U)\in\acc{w,w^R}$.
For example, avoiding $xyxy$ and $xyx^Ry$ simultaneously is equivalent to avoid $xyx^Uy$.
Notice that the notion of undirected avoidability recently considered by Currie and Mol~\cite{CM2020} corresponds to the case
where every occurrence of every variable of the pattern/formula is equipped with $-^U$.

Recall that a word is $(\beta^+,n)$-free if it contains no repetition with exponent strictly greater than $\beta$ and period at least $n$.
Also, a word is $(\beta^+)$-free if it is $(\beta^+,1)$-free.

\section{Formulas with unbounded avoidability index}\label{pphi}
Let us first show that for every $k\ge2$, $\phi_k$ is avoided by the periodic word $(\ell_0\ell_1\ldots\ell_k)^\omega$ over $(k+1)$ letters.
This word is 2-directed, so every occurrence $h$ of $\phi_k$ is such that $|h(x_i)|=1$ for every $0\le i<k$ by Remark~\ref{remark}.
Without loss of generality, $h(x_0)=\ell_0$. This forces $h(x_1)=\ell_1$, $h(x_2)=\ell_2$, and so on until
$h(x_{k-1})=\ell_{k-1}$ and $h(x_0)=\ell_k$, which contradicts $h(x_0)=\ell_0$.
Thus $\lambda(\phi_k)\le k+1$.

Let $b$ be an integer and let $w$ be an infinite word on at most $b$ letters.
Consider the Rauzy graph $R$ of $w$ such that the vertices of $R$ are the letters of $w$
and for every factor $uv$ of length two in $w$, we put the arc $\overrightarrow{uv}$ in $R$.
So $R$ is a directed graph, possibly with loops (circuits of length 1) and digons (circuits of length 2).
Since $w$ is infinite, every vertex of $R$ has out-degree at least 1. So $R$ contains a circuit $C_i$ of length $i$ with $1\le i\le b$.
Let $c_0,c_1,\ldots,c_{i-1}$ be the vertices of $C_i$ in cyclic order. Let $k$ be the least common multiple of $1,2,\ldots,b$.
Since $i$ divides $k$, $w$ contains the occurrence $h$ of $\phi_k$ such that $h(x_j)=c_{j\pmod i}$ for every $0\le j<k$. Thus $\lambda(\phi_k)>b$.

\section{Formulas that flatten to $xyzyx\cdot zyxyz$}\label{pxyzyx}
Notice that avoiding simultaneously the formulas in Theorem~\ref{xyzyx} is equivalent to avoiding $F=xyzy^Ux\cdot zy^Uxy^Uz\cdot y^R$.
The fragment $y^R$ is here to exclude the classical formula $xyzyx\cdot zyxyz$. Indeed, even if $xyzyx\cdot zyxyz$ is known to be 2-avoidable~\cite{circular},
a computer check shows that $xyzyx\cdot zyxyz$ and $F$ cannot be avoided simultaneously over two letters,
that is, $xyzy^Ux\cdot zy^Uxy^Uz$ is not 2-avoidable.


We use the method in~\cite{Ochem2004} to show that the image of every $(7/4^+)$-free word over $\Sigma_4$ by the following $21$-uniform morphism is $\paren{\tfrac{22}{15}^+,85}$-free.
We also check that such a binary word is $11$-directed.
$$
\begin{array}{ll}
\texttt{0}\to&\texttt{000010111000111100111}\\
\texttt{1}\to&\texttt{000010110011011110011}\\
\texttt{2}\to&\texttt{000010110001111010011}\\
\texttt{3}\to&\texttt{000010110001001101111}\\
\end{array}
$$
Consider an occurrence $h$ of $F$. Since $F$ contains $y\cdot y^R$,
then $|h(y)|\le10$ by Remark~\ref{remark}. Suppose that $|h(xz)|\ge83$. Then $h(xyzy^Ux)$ is a repetition with period $|h(xyzy)|\ge85$.
This implies $\frac{|h(xyzyx)|}{|h(xyzy)|}\le\frac{22}{15}$, which gives $|h(x)|\le\tfrac78|h(yzy)|$. Since $|h(y)|\le10$, we deduce $|h(x)|\le\tfrac{35}2+\tfrac78|h(z)|$.
Symmetrically, considering the repetition $h(zy^Uxy^Uz)$ gives $|h(z)|\le\tfrac{35}2+\tfrac78|h(x)|$.
So $$|h(x)|\le\tfrac{35}2+\tfrac78|h(z)|\le\tfrac{35}2+\tfrac78\paren{\tfrac{35}2+\tfrac78|h(x)|}=\tfrac{525}{16}+\tfrac{49}{64}|h(x)|$$
and $$|h(x)|\le\frac{\tfrac{525}{16}}{1-\tfrac{49}{64}}=140.$$
Symmetrically, $|h(z)|\le140$.

In every case, $|h(x)|\le140$, $|h(z)|\le140$, and $|h(y)|\le10$. Thus we can check exhaustively that $h$ does not exist.

\section{The formulas $xyzx\cdot yzxy\cdot z^R$ and $xyzx\cdot yz^Rxy$}\label{pxyzx}

Notice that avoiding $xyzx\cdot yzxy\cdot z^R$ and $xyzx\cdot yz^Rxy$ simultaneously is equivalent to avoiding $F=xyzx\cdot yz^Uxy\cdot z^R$.
We use the method in~\cite{Ochem2004} to show that the image of every $(7/4^+)$-free word over $\Sigma_4$ by the following $9$-uniform morphism is $\paren{\tfrac{131}{90}^+,28}$-free.
We also check that such a ternary word is $4$-directed.
$$
\begin{array}{ll}
\texttt{0}\to&\texttt{011122202}\\
\texttt{1}\to&\texttt{010121202}\\
\texttt{2}\to&\texttt{001112122}\\
\texttt{3}\to&\texttt{000101120}\\
\end{array}
$$
Consider an occurrence $h$ of $F$. Since $F$ contains $z\cdot z^R$,
then $|h(z)|\le3$ by Remark~\ref{remark}. Suppose that $|h(xy)|\ge27$. Then $h(xyzx)$ is a repetition with period $|h(xyz)|\ge28$.
This implies $\frac{|h(xyzx)|}{|h(xyz)|}\le\frac{131}{90}$, which gives $|h(x)|\le\tfrac{41}{49}|h(yz)|$. Since $|h(z)|\le3$, we deduce $|h(x)|\le\tfrac{129}{49}+\tfrac{41}{49}|h(y)|$.
Symmetrically, considering the repetition $h(yz^Uxyz)$ gives $|h(y)|\le\tfrac{129}{49}+\tfrac{41}{49}|h(x)|$.
So $$|h(x)|\le\tfrac{129}{49}+\tfrac{41}{49}|h(y)|\le\tfrac{129}{49}+\tfrac{41}{49}\paren{\tfrac{129}{49}+\tfrac{41}{49}|h(x)|}=\tfrac{11610}{2401}+\tfrac{1681}{2401}|h(x)|$$
and $$|h(x)|\le\frac{\tfrac{11610}{2401}}{1-\tfrac{1681}{2401}}=\tfrac{129}8.$$
So $|h(x)|\le16$ and, symmetrically, $|h(y)|\le16$.

In every case, $|h(xy)|\le32$ and $|h(z)|\le3$. Thus we can check exhaustively that $h$ does not exist.


\end{document}